\documentclass[a4paper,leqno]{article}

\usepackage{amsmath}
\usepackage{amsfonts}
\usepackage{enumerate}
\usepackage{theorem}

\theoremstyle{plain}
\newtheorem{teo}{Theorem}[section]
\newtheorem{lem}[teo]{Lemma}

\newtheorem{prop}[teo]{Proposition}

{\theorembodyfont{\rmfamily}
\newtheorem{defin}[teo]{Definition}

}

\renewcommand{\eqref}[1]{\textnormal{(\ref{#1})}}

\numberwithin{equation}{section}

\newcommand{\cvd}{\hfill$\square$}
\newcommand{\proof}[1]{\noindent\textsc{Proof#1}}
\newcommand{\rmi}{\mathrm{i}}
\newcommand{\rme}{\mathrm{e}}

\title{Determining a sound-soft polyhedral scatterer
 by a single far-field measurement\footnotemark[1]}

\author{\begin{tabular}{lcl}
Giovanni Alessandrini & and & Luca Rondi\\
\normalsize{\texttt{alessang@univ.trieste.it}} &
& \normalsize{\texttt{rondi@univ.trieste.it}}
\end{tabular}\\
\normalsize{Dipartimento di Scienze Matematiche}\\
\normalsize{Universit\`a degli Studi di Trieste, Italy}}

\date{}

\begin{document}

\maketitle
\footnotetext[1]{Work supported by MIUR under grant n.~2002013279.}

\setcounter{section}{0}
\setcounter{secnumdepth}{1}

\begin{abstract}
We prove that a sound-soft polyhedral scatterer is uniquely
determined by the far-field pattern corresponding to an incident
plane wave at one given wavenumber and one given incident
direction.
\end{abstract}

\bigskip

{\footnotesize \hspace{6cm}\begin{tabular}{l} Lo duca e io per
quel cammino ascoso\\intrammo a
ritornar nel chiaro mondo;\\ e sanza cura aver d'alcun riposo,\\
salimmo su, el primo e io secondo,\\ tanto ch'i' vidi de le cose
belle\\che porta'l ciel, per un pertugio tondo;\\ e quindi uscimmo
a riveder le stelle.\\ \\Dante, Inferno, C.XXXIV, 133-139.
\end{tabular}}

\section{Introduction}
We consider the acoustic scattering problem with a sound-soft
obstacle $D$. For simplicity of exposition, let us assume here
that $D$ is a bounded solid in $\mathbb{R}^N$, $N \geq 2$, that is
that $D$ is a connected compact set which coincides with the
closure of its interior.
We shall denote by $G$ the exterior of D
\begin{equation}
G=\mathbb{R}^N\backslash D
\end{equation}
and we shall assume throughout that it is connected.

Let $\omega\in \mathbb{S}^{N-1}$ and $k>0$ be fixed. Let $u$ be
the complex valued solution to
\begin{equation}\label{dirpbm}
\left\{\begin{array}{ll}
\Delta u + k^2u=0 & \text{in }G,\\
u(x)=u^s(x)+\rme^{\rmi k\omega\cdot x}& x\in G,\\
u=0 & \text{on }\partial G,\\
\lim_{r\to \infty}r^{(N-1)/2}\left(\frac{\partial u^s}{\partial
r}-\rmi ku^s\right)=0 & r=\|x\|.
\end{array}\right.
\end{equation}
It is well-known that the asymptotic behavior at infinity of
the so-called scattered field $u^s(x)=u(x)-\rme^{\rmi k\omega\cdot
x}$ is governed by the formula
\begin{equation}\label{asympt}
u^s(x)=\frac{\rme^{\rmi k\|x\|}}{\|x\|^{(N-1)/2}}\left\{
u_{\infty}(\hat{x})+O\left(\frac{1}{\|x\|}\right) \right\},
\end{equation}
as $\|x\|$ goes to $\infty$, uniformly in all directions
$\hat{x}=x/\|x\|\in \mathbb{S}^{N-1}$. The function $u_{\infty}$,
which is defined on $\mathbb{S}^{N-1}$, is called the
\emph{far-field pattern} of $u^s$, see for instance \cite{Col e
Kre98}. In this paper we prove that if $N=2$ and $D$ is a polygon,
or if $N=3$ and $D$ is a polyhedron, then it is uniquely
determined by the far-field pattern $u_{\infty}$ for one
wavenumber $k$ and one incident direction $\omega$, see Theorem
\ref{mainthm} below. Indeed we prove this result independently of
the dimension $N\geq 2$, and for this reason it is convenient to
express our assumption on $D$ by prescribing that it is an
$N$-dimensional polyhedron, that is, a solid whose boundary is
contained into the union of finitely many $(N-1)$-dimensional
hyperplanes (more precisely we should say a \emph{polytope}, see
for instance \cite{Cox}, but for the sake of simplicity we prefer
to stick to
the $3$-dimensional terminology). 
In fact our result applies to a wider class of scatterers D, which
need not to be solids, nor connected, but whose boundary is the
finite union of the closures of open subsets of
$(N-1)$-dimensional hyperplanes. See Section \ref{thmsec} below
for a complete formulation.

We wish to mention here that in '94  C. Liu and A. Nachman \cite
{Liu e Nach} proved, among various results, that, for $N\geq 2$,
$u_{\infty}$ uniquely determines the convex envelope of a
polyhedral obstacle $D$, and they also outlined a proof of the
unique determination of a polyhedral obstacle. Their arguments
involve a scattering theory analogue of a classical theorem of
Polya on entire functions and the reflection principle for
solutions of the Helmholtz equation across a flat boundary.

More recently J. Cheng and M. Yamamoto \cite{Che e Yam}, for the
case $N=2$, proved that the far-field pattern uniquely determines
a polygonal obstacle $D$, provided $D$ satisfies an additional
geometrical condition, which, roughly speaking, is expressed in
terms of the absence of trapped rays in its exterior $G$. The
method of proof in \cite{Che e Yam} is mainly based on the use of
the reflection principle and on the study of the behavior of the
nodal line $\{u=0\}$ of the solution $u$ to (\ref{dirpbm}) near the
boundary $\partial G$.

Also in this paper we make use of such a reflection argument, but,
rather than examining the boundary behavior of the nodal set we
investigate the structure of the nodal set of $u$ in the interior
of $G$. In this respect, the main tool is summarized in the fact
that if $D$ is a polyhedron, then the nodal set of $u$ in $G$ does
not contain any open portion of an $(N-1)$-dimensional hyperplane,
see Theorem \ref{noflatpointthm}.

In the next Section \ref{thmsec} we set up our main hypotheses on
the obstacle, we state the main results Theorem \ref{mainthm} and
Theorem \ref{noflatpointthm} and prove Theorem \ref{mainthm}.

In Section \ref{pathsec} we prove Theorem \ref{noflatpointthm}.
The proof is preceded by a sequence of Propositions and auxiliary
Lemmas regarding the study of the nodal sets of real valued
solutions to the Helmholtz equation,  see Proposition
\ref{orderprop}, and the construction of a suitable path in $G$
(\emph{cammino ascoso} $=$ \emph{hidden path}) which connects a point in
$\partial D$ to infinity, avoiding the singular points in the
nodal set of $u$ and intersecting the nodal set orthogonally,
Proposition \ref{complexpathprop}.

\section{The uniqueness result}\label{thmsec}
\begin{defin}
Let us define a \emph{cell} as the closure of an open subset of an
$(N-1)$-dimensional hyperplane. We shall say that $D$ is a
\emph{polyhedral scatterer}
  if it is a compact subset of $\mathbb{R}^N$, such that
\begin{enumerate}[(i)]
\item the exterior $G=\mathbb{R}^N\backslash D$ is connected,
\item the boundary of $G$ is given by the finite union of cells
$C_j$. \label{cell}
\end{enumerate}
\end{defin}

Let us observe that an equivalent condition to (\ref{cell}) is
that $D$ has the form
\begin{equation*}
D=(\bigcup_{i=1}^mP_i)\cup(\bigcup_{j=1}^nS_j),
\end{equation*}
where each $P_i$ is a polyhedron and each $S_j$ is a cell, thus we
are allowing the simultaneous presence of solid obstacles and of
crack-type scatterers. Note also that, by this definition, a cell
needs not to be an $(N-1)$-dimensional polyhedron.

We also recall that for any compact set $D$ a weak solution $u \in
W^{1,2}_{loc}(G)$ to (\ref{dirpbm}) exists and is unique, see for
instance \cite{Ram e Rui}. As is well-known, $u$ is analytic in
$G$, but, of course, due to the possible irregularity of the
boundary of $G$, the Dirichlet boundary condition in
(\ref{dirpbm}) is, in general, satisfied in the weak sense only.
On the other hand, one can notice that, if $x^0 \in
\partial G$ is an interior point of one of the cells forming
$\partial G$, then it is a regular point for the Dirichlet problem
in $G$, hence $u$ is continuous up to $x^0$ and $u(x^0)=0$.

\begin{teo}\label{mainthm}
Let us fix $\omega\in \mathbb{S}^{N-1}$ and $k>0$. A
polyhedral scatterer $D$ is uniquely determined by the far-field
pattern $u_{\infty}$.
\end{teo}

A proof of Theorem~\ref{mainthm} will be obtained as a consequence
of Theorem \ref{noflatpointthm} below, the following definitions
will be needed.

\begin{defin}
Let us denote by $\mathcal{N}_u$ the \emph{nodal set} of $u$ in
$G$, that is
\begin{equation*}
\mathcal{N}_u=\{x\in G:\ u(x)=0\}.
\end{equation*}
We shall say that $x\in \mathcal{N}_u$ is a \emph{flat point} if
there exist a hyperplane $\Pi$ through $x$ and a positive number
$r$ such that $\Pi\cap B_r(x)\subset \mathcal{N}_u$.
\end{defin}

\begin{teo}\label{noflatpointthm}
Let $D$ be a polyhedral scatterer. Then $\mathcal{N}_u$ cannot
contain any flat point.
\end{teo}

We postpone the proof of this result to the next
Section~\ref{pathsec} and we conclude the proof of
Theorem~\ref{mainthm}.

\bigskip

\proof{ of Theorem~\ref{mainthm}.} Let $D$ and $D'$ be two
polyhedral scatterers and let $u'$ be the solution to (\ref{dirpbm})
when $D$ is replaced with $D'$. Let us assume that for a given
$\omega\in \mathbb{S}^{N-1}$ and $k>0$, $u_{\infty}=u'_{\infty}$.
We denote with $\widetilde{G}$ the connected component of
$\mathbb{R}^N\backslash(D\cup D')$ which contains the exterior of
a sufficiently large ball. By Rellich's Lemma (see for instance
\cite[Lemma~2.11]{Col e Kre98}) and unique continuation we infer
that $u=u'$ over $\widetilde{G}$.

First, we notice that if $\partial\widetilde{G}\subset D\cap D'$,
then $D=D'=\mathbb{R}^N\backslash\widetilde{G}$. This is due to
the fact that both $G$ and $G'=\mathbb{R}^N\backslash D'$ are
connected.

Let us assume, by contradiction, that $D$ is different from $D'$.
Then, without loss of generality, we can assume that there exists
a point $x'\in(\partial G'\backslash D) \cap
\partial\widetilde{G}$. We can also assume that $x'$ belongs to
the interior of one of the cells composing $\partial G'$, and
therefore that there exist a hyperplane $\Pi'$ and $r>0$ such that
$x'\in S'=\Pi' \cap B_r(x')\subset (\partial G'\backslash D) \cap
\partial\widetilde{G}$. Since $u=u'$ in $\widetilde{G}$, by continuity we have that $u=u'=0$ on $S'$, hence $S'$ is
contained into the nodal set of $u$, that is $S'\subset
\mathcal{N}_u$, and, consequently, $x'$ is a flat point for
$\mathcal{N}_u$. This contradicts
Theorem~\ref{noflatpointthm}.\cvd

\section{The \emph{hidden path} and the proof
of Theorem~\ref{noflatpointthm}}\label{pathsec}

We start with a well-known property of the nodal set of $u$.

\begin{lem}\label{Nbounded}
The nodal set $\mathcal{N}_u$ is bounded.
\end{lem}

\proof{.} By \eqref{asympt}, we have that the scattered field
$u^s(x)$ tends to zero, as $\|x\|$ tends to infinity, uniformly
for all directions $\hat{x}=x/\|x\|\in \mathbb{S}^{N-1}$. Then the
lemma immediately follows by observing that
$|u(x)|=|u^s(x)+\rme^{\rmi k\omega\cdot x}| \rightarrow 1$
uniformly as $\|x\| \rightarrow \infty$.\cvd

\bigskip

Next we discuss some properties of the nodal set of real valued
solutions to the Helmholtz equation. Let $v$ be a nontrivial real
valued solution to
\begin{equation}\label{eq}
\Delta v + k^2v=0 \ \text{in }G ,
\end{equation}
in a connected open set $G$.
We denote the \emph{nodal set} of $v$
as
$$\mathcal{N}_v=\{x\in G:\ v(x)=0\}$$
and we let $\mathcal{C}_v$ be the set of \emph{nodal critical
points}, that is
$$\mathcal{C}_v=\{x\in G:\ v(x)=0\text{ and }\nabla v(x)=0\}.$$
We say that $\Sigma\subset \mathcal{N}_v$ is a \emph{regular
portion} of $\mathcal{N}_v$ if it is an analytic open and
connected hypersurface contained in
$\mathcal{N}_v\backslash\mathcal{C}_v$. Let us denote by
$A_1,A_2,\ldots,A_n,\ldots$ the \emph{nodal domains} of $v$,
that is the connected components of $\{x\in G:\ v(x)\neq
0\}=G\backslash \mathcal{N}_v$.

\begin{prop}\label{orderprop}
We can order the nodal domains $A_1,A_2,\ldots,A_n,\ldots$ in such
a way that for any $j\geq 2$ there exist $i$, $1\leq i<j$, and a
regular portion $\Sigma_j$ of $\mathcal{N}_v$ such that
\begin{equation}\label{order}
\Sigma_j\subset \partial A_i\cap\partial A_j.
\end{equation}
\end{prop}

We subdivide the main steps of the proof of this proposition in
the next two lemmas.

\begin{lem}\label{lemma1}
Let $A_1,\ldots,A_n$ be nodal domains and let
$A=\stackrel{\circ}{\overline{A_1\cup\ldots\cup A_n}}$. If $x\in
\partial A\cap G$, then for any $r>0$ there exists $y\in
(B_r(x)\cap G)\backslash\overline{A}$.
\end{lem}

\proof{.} We can assume, without loss of generality, that $r>0$ is
such that $B_r(x)\subset G$. Then, let us assume, by
contradiction, that we have $B_r(x)\subset\overline{A}$. Then we
infer that $x\in\stackrel{\circ}{\overline{A}}=A$ and this
contradicts the fact that $x\in\partial A$.\cvd

\begin{lem}\label{lemma2}
Let $A_1,\ldots,A_n$ be nodal domains and let
$A=\stackrel{\circ}{\overline{A_1\cup\ldots\cup A_n}}$. If $x\in
\partial A\cap G$, then for any $r>0$ there exists $y\in
B_r(x)\cap \partial A\cap G$ such that $\nabla v(y)\neq 0$.
\end{lem}

\proof{.} We can assume, without loss of generality, that $r>0$ is
such that $B_r(x)\subset G$. Assume, by contradiction, that
$\nabla v\equiv 0$ on $B_r(x)\cap \partial A$ and set $w=v$ in
$B_r(x)\cap A$, $w=0$ in $B_r(x)\backslash A$. One can easily
verify that $w\in W^{2,\infty}(B_r(x))$ and also that $w$ is a
strong solution to the Helmholtz equation in $B_r(x)$. Now, by
Lemma~\ref{lemma1}, $w\equiv 0$ on an open subset of $B_r(x)$ and
hence by unique continuation $w\equiv 0$ in $B_r(x)$ which is
impossible.\cvd

\bigskip

\proof{ of Proposition~\ref{orderprop}.}
We proceed by induction. We choose $A_1$ arbitrarily.

Let us assume that we have ordered $A_1,\ldots,A_n$ in such a way
that there exist $\Sigma_2,\ldots,\Sigma_n$ regular portions of
$\mathcal{N}_v$ such that \eqref{order} holds for any
$j=2,\ldots,n$ and for some $i<j$.

Let $A=\stackrel{\circ}{\overline{A_1\cup\ldots\cup A_n}}$. If
$A=G$, then we are done. Otherwise, since $G$ is connected, we can
find $x\in\partial A\cap G$. We apply Lemma~\ref{lemma2} and we
fix, for $r>0$ small enough, a point $y\in B_r(x)\cap \partial
A\cap G$ such that $\nabla v(y)\neq 0$. There exists a positive
$r_1$ such that $B_{r_1}(y)\cap \partial A$ is a regular portion
of $\mathcal{N}_v$ and there exist exactly two nodal domains,
$\tilde{A}_1\subset A$ and $\tilde{A}_2$ with $\tilde{A}_2\cap
A=\emptyset$, whose intersections with $B_{r_1}(y)$ are not empty.
It is clear that $\tilde{A}_1$ coincides with $A_i$ for some
$i=1,\ldots,n$, and if we pick $A_{n+1}=\tilde{A}_2$ and choose
$\Sigma_{n+1}=B_{r_1}(y)\cap \mathcal{N}_v$, then \eqref{order}
holds for $j=n+1$, too.\cvd

\bigskip

We now show that we are able to connect points of
$G\backslash\mathcal{N}_v$ with suitable regular curves contained
in $G$ which avoid the nodal critical points of $v$. Here and in the sequel we shall say
that a curve $\gamma=\gamma(t)$ is regular if it is $C^1$-smooth and
$\frac{\mathrm{d}}{\mathrm{d}t}\gamma(t)\neq 0$
for every $t$.

\begin{prop}\label{realpathproof}
Let $x_1$ and $y_1$ belong to $G\backslash \mathcal{N}_v$. Then
there exists a regular curve $\gamma$ contained in $G$ and
connecting $x_1$ with $y_1$ such that the following conditions are
satisfied
\begin{equation}\label{nocritic}
\gamma\cap\mathcal{C}_v = \emptyset,
\end {equation}
 \begin{equation}\label{orthogon}\text{if }x \in
\mathcal{N}_v\cap\gamma , \text{then } \gamma \text{ intersects }
\mathcal{N}_v  \text{ at }x  \text{ orthogonally.}
\end{equation}
\end{prop}

\proof{.} We order the nodal domains $A_1,A_2,\ldots,A_n,\ldots$ according
to Proposition~\ref{orderprop}. Without loss of
generality we can assume that $x_1\in A_1$ and $y_1\in A_i$ for
some $i>1$. By Proposition~\ref{orderprop}, we can find $i_l$,
with $l=1,\ldots,n$, such that $i_1=1$, $i_n=i$, and, for any
$l=2,\ldots,n$, $i_{l-1}<i_l$ and there exists a regular portion
of $\mathcal{N}_v$, $\Sigma_{i_l}$, such that
$\Sigma_{i_l}\subset\partial A_{i_{l-1}}\cap\partial A_{i_l}$.

Let $\sigma_l$ be a line segment crossing $\Sigma_{i_l}$
orthogonally and let it be small enough such that $\sigma_l
\subset A_{i_{l-1}}\cup\Sigma_{i_l}\cup A_{i_l}$. Let $y_l^{-}\in
A_{i_{l-1}}, \ y_l^{+}\in A_{i_l}$ be the endpoints of $\sigma_l$.
Let $\beta_1$ be a regular path within $A_1$ which joins $x_1$ to
$y_2^{-}$ and has a
$C^1$-smooth junction with $\sigma_2$ at $y_2^{-}$.
For every $l=2,\ldots,n-1$, let $\beta_l$ be a
regular path within $A_{i_l}$ which joins $y_l^{+}$ to
$y_{l+1}^{-}$ and has $C^1$-smooth junctions with the segments
$\sigma_l$ and $\sigma_{l+1}$, at the points $y_l^{+}$,
$y_{l+1}^{-}$, respectively. Let $\beta_n$ be a regular
path within $A_{i_n}$ which joins $y_n^{+}$ to $y_1$ and has a
$C^1$-smooth junction with $\sigma_n$ at $y_n^{+}$. We form the curve
$\gamma$ by attaching consecutively the curves $\beta_1, \
\sigma_1, \ \beta_2, \ \sigma_2, \ldots$ up to $\beta_n$.\cvd

\bigskip

We have what is needed to build up our \emph{hidden path}. From
now on we consider $G=\mathbb{R}^N\backslash D$ and $v=\Re u$.
Note that $\mathcal{N}_u \subset \mathcal{N}_v$.

\begin{prop}\label{complexpathprop}
Let $x_1\in \partial G$ be such that $x_1$ belongs to the interior
of one of the cells forming $\partial G$ and $\frac{\partial
v}{\partial \nu}(x_1)\neq 0$, $\nu$ being the unit normal to
$\partial G$ at $x_1$, pointing to the interior of $G$. Let
$y_1\in \mathcal{N}_u \backslash \mathcal{C}_v$ be fixed.

Then there exists a regular curve
$\gamma:[0,+\infty)\mapsto\mathbb{R}^N$, such that the following
conditions are satisfied

\begin{enumerate}[1\textnormal{)}]

\item $\gamma(0)=x_1$\textnormal{;}

\item $\gamma(t)\in G$ for every $t>0$\textnormal{;}

\item there exists $t_1$ such that $\gamma(t_1)=y_1$\textnormal{;}

\item $\lim_{t\to+\infty}\|\gamma(t)\|=+\infty$\textnormal{;}

\item if, for some $t$, $\gamma(t)\in \mathcal{N}_u$,
then $\gamma(t)\notin \mathcal{C}_v$  and $\gamma$ intersects
$\mathcal{N}_v$ at $x=\gamma(t)$ orthogonally.

\end{enumerate}

\end{prop}

\proof{.} Let $A_1$ be the nodal domain of $v$ such that $x_1 \in
\partial A_1$ and let $\eta_1$ be a line segment in $A_1$ having
$x_1$ as an endpoint and which is orthogonal to $\partial G$
there. Let $x'_1\in A_1$ be the other endpoint of $\eta_1$. Let
$\eta_2$ be a line segment crossing $\mathcal{N}_v$ orthogonally
at the point $y_1$. Let it be small enough so that $v$ is strictly
monotone on $\eta_2$. Let $y'_1,\ y''_1$ be the endpoints of
$\eta_2$. By Proposition~\ref{realpathproof}, we can find a
regular curve $\gamma_1$ joining $x'_1$ to $y'_1$ and satisfying
conditions (\ref{nocritic}), (\ref{orthogon}). We can also choose
$\gamma_1$ in such a way that it has $C^1$-smooth junctions with the
segments $\eta_1, \ \eta_2$ at its endpoints. Let $R>0$ be large
enough so that $\mathcal{N}_u \subset B_R(0)$ and let us fix $z_1,
\ |z_1|>R$.  Again by Proposition~\ref{realpathproof}, we can find
a regular curve $\gamma_2$ joining $y''_1$ to $z_1$ and
satisfying conditions (\ref{nocritic}), (\ref{orthogon}) and also
such that it has a $C^1$-smooth junction with $\eta_2$ at the point
$y''_1$. Next let us fix a regular path $\gamma_3$ in
$\mathbb{R}^N\backslash B_R(0)$ joining $z_1$ to $\infty$ having a
$C^1$-smooth junction with $\gamma_2$ at $z_1$. The resulting path
$\gamma$ is obtained by attaching the paths $\eta_1, \ \gamma_1, \
\eta_2, \ \gamma_2, \ \gamma_3$. \cvd

\begin{lem}\label{nextflatpointlem}
Let the assumptions of Proposition~\textnormal{\ref{complexpathprop}} be
satisfied and let $\gamma$ be the path constructed there. If $y'=
\gamma(t')\in\mathcal{N}_u$ is a flat point, then there exists
$t''>t'$ such that $y''=\gamma(t'')\in \mathcal{N}_u$ is a flat
point.
\end{lem}
\proof{.} Let $\Pi'$ be the plane through $y'$ and let $r>0$ be
such that $S'=\Pi'\cap B_r(y')\subset \mathcal{N}_u$.

Let $\widetilde{S}'$ be the connected component of $\Pi'\backslash
D$ containing $y'$. We have that, by analytic continuation, $u$ is
identically zero on $\widetilde{S}'$. Therefore, we can
immediately notice that, by Lemma~\ref{Nbounded}, $\widetilde{S}'$
is bounded.

Let $\epsilon >0$ be small enough so that $v(\gamma(t))$ is
strictly monotone for $t'-\epsilon \leq t \leq t'+\epsilon$, and
let us set $y^- =\gamma(t'-\epsilon), \ y^+ =\gamma(t'+\epsilon)$.

Let $G^+$ be the connected component of
$G\backslash\widetilde{S}'$ containing $y^+$ and let $G^-$ be the
connected component of $G\backslash\widetilde{S}'$ containing
$y^-$. Let us remark that it may happen that $G^+=G^-$.

Let us denote with $R$ the reflection in $\Pi'$. We call $E^+$ the
connected component of $G^+\cap R(G^-)$ containing $y^+$ and $E^-$
the connected component of $G^-\cap R(G^+)$ containing $y^-$. We
observe that $E^-=R(E^+)$ and we set $E=E^+\cup E^-\cup
\widetilde{S}'$.

We have that $E$ is a connected open set and, by construction, the
boundary of $E$ is composed by cells, more precisely by subsets of
the cells of $\partial G$ and of $R (\partial G)$. Furthermore, in
$E$ we have that $u=-Ru$ where $Ru(x)=u(R(x))$. In fact, $u+Ru$ is
a solution of the Helmholtz equation  in $E$ with zero Cauchy data
on $\widetilde{S}'$.

In other words, $u$ is odd symmetric in $E$, with respect to the
plane $\Pi'$. Hence, we infer that $u=0$ on $E\cap\Pi'$
and, moreover, $u$ is continuous up to the interior of each cell
forming  $\partial E$ and $u=0$ there. Furthermore, since $u$ is
continuous in $G$, we have that $u=0$ in all of $\partial E \cap
G$. That is $\partial E \cap G \subset \mathcal{N}_u$.

Let us exclude now the case that $E$ is unbounded. In fact,
$\partial E$ is bounded and, if $E$ were unbounded, then $E$ would
contain $\mathbb{R}^N\backslash B_{\rho}(0)$, for some
sufficiently large $\rho>0$. Then $u=0$ on $\Pi'\backslash
B_{\rho}(0)$ and this contradicts Lemma~\ref{Nbounded}.

Thus $E$ is a bounded open set containing $y'$. Since $\gamma$ is
not bounded, there exists $t''>t'$ such that
$\gamma(t'')\in\partial E\cap G$. We have that
$y''=\gamma(t'')\in\mathcal{N}_u$ and, by the properties of
$\gamma$, it is not a critical point of $v$. Let $C$ be a cell of
$\partial E$ such that $y'' \in C$ and let $\Pi''$ be the
hyperplane containing $C$. Let $r>0$ be such that $B_r(y'')
\subset G$. We have that $u=0$ on $C \cap B_r(y'')$ and hence, by
analytic continuation, $u=0$ on $\Pi'' \cap B_r(y'')$, therefore
$\Pi'' \cap B_r(y'') \subset \mathcal{N}_u $.\cvd

\bigskip

\proof{ of Theorem~\ref{noflatpointthm}.} Let us assume, by
contradiction, that $y_1\in \mathcal{N}_u$ is a flat point. Let
$\Pi_1$ be the plane through $y_1$ and $r>0$ such that
$S_1=\Pi_1\cap B_r(y_1)\subset \mathcal{N}_u$. By the uniqueness
for the Cauchy problem, $S_1$ contains at least one point
$y_1'\notin \mathcal{C}_v$. Thus, without loss of generality, we
can assume that there exists a flat point $y_1\in \mathcal{N}_u
\backslash  \mathcal{C}_v$.

We arbitrarily fix a point $x_1$ belonging to the interior of one
of the cells of $\partial G$. Again by the uniqueness
for the Cauchy problem, we can
assume, without loss of generality, that $\frac{\partial
v}{\partial \nu}\neq 0$, $\nu$ being the interior unit normal to
$\partial G$ at the point $x_1$.

We choose $\gamma$  according to
Proposition~\ref{complexpathprop}. Then, applying iteratively
Lemma \ref{nextflatpointlem}, we can find a strictly increasing
sequence $\{t_n\}_{n\in\mathbb{N}}$ such that, for any $n$,
$y_n=\gamma(t_n)$ is a flat point of $u$ and, by construction of
$\gamma$, $y_n$ is not a critical point of $v$. Since
$\mathcal{N}_u$ is bounded and $\lim_{t\to
+\infty}\|\gamma(t)\|=+\infty$, there exists a finite $T$ such
that $\lim_{n\to +\infty}t_n=T$. We have that
$\tilde{y}=\gamma(T)$ belongs to $\mathcal{N}_u$ and, again by the
properties of $\gamma$, $\tilde{y}$ is not a critical point of $v$
and $\gamma$ is orthogonal to $\mathcal{N}_v$ there. Therefore,
there exists $\delta >0$ such that $v(\gamma(t)) \neq 0$ for every
$ T-\delta <t < T$ and this contradicts the fact that
$\gamma(t_n)\in\mathcal{N}_u$ for any $n$.\cvd

\subsection{Acknowledgements}
The authors wish to express their gratitude
to Adrian Nachman and Masahiro Yamamoto for kindly sending
their respective preprints \cite{Liu e Nach} and \cite{Che e Yam}. The research reported
in this paper originated at the 2003 Oberwolfach meeting 
``Inverse Problems in Wave Scattering and Impedance Tomography'', the authors wish to thank the organizers
Martin Hanke-Bourgeois,
Andreas Kirsch and William Rundell,
and the Mathematisches Forschungsinstitut Oberwolfach.

\end{document}